\documentclass[sts]{imsart}

\RequirePackage{amsthm,amsmath,amsfonts,amssymb}
\RequirePackage[authoryear]{natbib}
\RequirePackage[colorlinks,citecolor=blue,urlcolor=blue,linkcolor=blue]{hyperref}

\startlocaldefs

\numberwithin{equation}{section}
\allowdisplaybreaks
\AtBeginDocument{\sloppy}

\newcommand{\E}{\mathbb{E}}
\newcommand{\N}{\mathcal{N}}
\newcommand{\D}{\mathcal{D}}
\newcommand{\I}{\mathcal{I}}
\DeclareMathOperator{\tr}{tr}
\DeclareMathOperator{\KL}{KL}
\DeclareMathOperator{\Var}{Var}
\DeclareMathOperator{\argmax}{arg\,max}
\DeclareMathOperator{\argmin}{arg\,min}

\endlocaldefs

\begin{document}

\begin{frontmatter}

\title{A Bayesian-Geometric Reading of Classical Model-Based Optimal Design}
\runtitle{Bayesian-Geometric Reading of Optimal Design}

\begin{aug}
\author{\fnms{Mohamed}~\snm{Tarek}\ead[label=e1]{mohamed@pumas.ai} \orcid{0000-0003-3074-0678}}
\address{PumasAI, USA\printead[presep={\ }]{e1}.}
\end{aug}

\begin{abstract}
The classical model-based optimal design criteria (D-, T-, L-, C-, and A-optimality) are usually motivated via the Cram\'er--Rao bound as optimizing scalar functions of the Fisher information matrix (FIM). This paper offers a complementary reading in which each criterion appears asymptotically as $N \to \infty$ as the leading-order term of an expected dissimilarity between prior and posterior, with the FIM entering through the Bernstein--von Mises limit. The Kullback--Leibler (KL) divergence of the posterior relative to the prior recovers D-optimality (classically known); a reverse KL of the prior relative to the posterior recovers T-optimality (here defined as maximization of the trace of the FIM, following one common convention) and its weighted-trace generalization; and a first-order approximation of the squared 2-Wasserstein distance recovers L-, C- and A-optimality. The five criteria emerge as instances of a common Bayesian-geometric template, indexed by three choices: the dissimilarity, the prior parametrization, and the approximation scheme used to evaluate it (e.g.\ Taylor expansions of the matrix square root). For the weighted criteria, the Gaussian prior covariance encodes the criterion's weight matrix, so any choice of weights corresponds to a specific prior in the family. The T-, L-, C-, and A-optimality derivations appear novel to the best of my knowledge. A sixth criterion, E-optimality, is recovered not as an independent expected-dissimilarity limit but as a min-max robustification of C-optimality against uncertainty in the modeler's prior direction. The framework is extended to distinguish between a designer's prior over unknown parameters at design time and a modeler's prior used for posterior inference, recovering the classical expected information gain objective as a special case.
\end{abstract}

\begin{keyword}
\kwd{optimal experimental design}
\kwd{Fisher information}
\kwd{Bayesian experimental design}
\kwd{Kullback--Leibler divergence}
\kwd{Wasserstein distance}
\end{keyword}

\end{frontmatter}

\section{Introduction}
\label{sec:introduction}

Optimal design selects an experimental configuration $\delta$ from a feasible set $\D$ so as to maximize information gained about the parameters $\theta \in \mathbb{R}^k$ of a probabilistic model from observed data \citep{kiefer1959,fedorov1972,silvey1980,pukelsheim1993,atkinson2007,pronzatopazman2013}. Throughout this paper, $y$ denotes the random variable associated with a single observation under design $\delta$, with likelihood $p(y \mid \theta, \delta)$. A dataset of $N$ independent observations of $y$ under design $\delta$ is denoted by $D = \{y_i\}_{i=1}^N$, with joint likelihood $p(D \mid \theta, \delta) = \prod_{i=1}^N p(y_i \mid \theta, \delta)$. We write $\theta_0$ for the true (but unknown) value of the parameters that generated the data, $\hat{\theta}_N(D, \delta)$ for the maximum likelihood estimator of $\theta$ from $D$ under design $\delta$, $\I$ for the $k \times k$ identity matrix, and $\lambda_{\min}(M)$ (resp.\ $\lambda_{\max}(M)$) for the smallest (resp.\ largest) eigenvalue of a symmetric matrix $M$. The per-observation Fisher information matrix (FIM) is
\begin{align*}
F(\theta, \delta) &= \E_{y \mid \theta, \delta}\!\Bigl[
\Bigl(\frac{\partial}{\partial \theta} \log p(y \mid \theta, \delta)\Bigr) \\
&\qquad \Bigl(\frac{\partial}{\partial \theta} \log p(y \mid \theta, \delta)\Bigr)^{\!T}
\Bigr],
\end{align*}
and the prior, prior predictive, and posterior are denoted by $p(\theta)$, $p(D \mid \delta) = \int p(D \mid \theta, \delta)\,p(\theta)\,\mathrm{d}\theta$, and $p(\theta \mid D, \delta)$, respectively.

The classical model-based criteria (D-, T-, L-, C-, and A-optimality) are typically stated as optimizing scalar functions of $F(\theta_0, \delta)$, the locally optimal design setup of \citet{chernoff1953}:
\begin{align*}
\delta^*_D &= \argmax_{\delta \in \D} \log \det F(\theta_0, \delta) \\
\delta^*_T &= \argmax_{\delta \in \D} \tr(F(\theta_0, \delta)) \\
\delta^*_L &= \argmin_{\delta \in \D} \tr(L F(\theta_0, \delta)^{-1} L^T) \\
\delta^*_A &= \argmin_{\delta \in \D} \tr(F(\theta_0, \delta)^{-1}) \\
\delta^*_C &= \argmin_{\delta \in \D} c^T F(\theta_0, \delta)^{-1} c,
\end{align*}
where $L$ is a weight matrix and $c^T$ is a weight row vector; the $c^T F^{-1} c$ functional underlying C-optimality is due to \citet{elfving1952}, developed in the context of linear regression where $F$ does not depend on $\theta_0$. The standard motivation is the asymptotic normality of the maximum likelihood estimator \citep{vandervaart1998}, whose asymptotic covariance is $N^{-1} F(\theta_0, \delta)^{-1}$: each criterion is a particular way of shrinking this covariance. This paper offers a complementary, Bayesian-geometric reading, in the tradition of Bayesian experimental design \citep{lindley1956,degroot1962,chaloner1995,ryan2016,rainforth2024}. The model structure is assumed to be true, but the parameters are unknown. A single set of unknown parameters $\theta_0$ is assumed to have generated the data, and the goal is to learn about $\theta_0$ from $D$.

The modeler has a prior distribution $p(\theta)$ over the values of $\theta$, encoding beliefs about the parameters before observing any data. After the data is observed, the modeler updates their beliefs using Bayes' rule,
\[
p(\theta \mid D, \delta) = \frac{p(D \mid \theta, \delta)\,p(\theta)}{p(D \mid \delta)},
\]
to obtain the posterior $p(\theta \mid D, \delta)$. The design $\delta$ determines the likelihood $p(D \mid \theta, \delta)$ but the prior $p(\theta)$ is independent of the design. The prior predictive $p(D \mid \delta)$ is not generally the same as the true data-generating distribution of $D$ under design $\delta$, which is $p(D \mid \theta_0, \delta)$, because the prior may be misspecified.

In this setting, the experiment is informative to the extent that it has the potential to significantly change our beliefs about $\theta$ from the prior to the posterior. Let $d(p_0, p_1) \geq 0$ be a measure of dissimilarity on the space of distributions over $\theta$, with $d(p_0, p_1) = 0$ only when $p_0 = p_1$ almost everywhere; $d$ may be either a proper distance metric or an asymmetric divergence. The potential of an experiment to change our beliefs can then be quantified as the expected dissimilarity between the prior and the posterior, where the expectation is taken over the true data-generating distribution of $D$ under design $\delta$
\[
U(\delta) = \E_{D \mid \theta_0, \delta}\bigl[d(p(\theta), p(\theta \mid D, \delta))\bigr].
\]
Since $D$ is unobserved at design time, one can ask which design would, on average over hypothetical repetitions of the experiment, generate the largest expected change in belief about $\theta$, i.e., the largest expected dissimilarity between the prior and the posterior. The optimal design is then
\[
\delta^* = \argmax_{\delta \in \D} U(\delta).
\]
The same construction applies to any derived quantity $q = g(\theta)$ of interest (e.g.\ a transformation of the parameters or a future observation): replacing $p(\theta)$ and $p(\theta \mid D, \delta)$ with the corresponding prior predictive distribution $p(q)$ and posterior predictive distribution $p(q \mid D, \delta)$ yields an analogous criterion targeting that quantity.

As a motivating example, consider selecting sample times in a pharmacokinetic experiment to estimate a drug's absorption and elimination parameters \citep{mentre1997}. Suppose the drug is given at time 0 and is fully cleared well before 24 hours. As a modeler, I may choose to have a prior over the absorption and elimination rates that predicts that concentrations fall below the limit of quantification at 48 hours. This prior is not too informative, in the sense that it is diffuse over a wide range of absorption and elimination rates, yet sufficiently informative to exclude sustained high levels after 48 hours. Given that prior, a non-informative design would measure concentrations 2 or more days after dosing. The data collected would merely confirm the near-zero drug concentration prediction of the prior (up to measurement noise). Therefore, the posterior would remain essentially identical to the prior. An informative design instead samples within the first 24 hours, when concentrations vary rapidly and different parameters produce markedly different trajectories. The data can then discriminate among a-priori plausible values, driving the posterior to concentrate on a much smaller region of the parameter space than the prior. The expected dissimilarity between prior and posterior over all possible datasets $D$ quantifies precisely this potential for an experiment to change our beliefs.

More broadly, this viewpoint admits a range of criteria based on choices of:
\begin{itemize}
\item the dissimilarity measure on distributions over $\theta$,
\item the prior distribution, and
\item the approximation scheme used to evaluate the dissimilarity, e.g.\ asymptotic limits or Taylor expansions.
\end{itemize}
Five classical (D-, T-, L-, C-, and A-) criteria arise as particular choices of these three ingredients, where the FIM enters through the Bernstein--von Mises (BvM) limit \citep{bernstein1917,vonmises1931,lecam1953,vandervaart1998,fischer2025} (formally introduced in Section~\ref{sec:setup}), exposing each criterion as a particular notion of geometric separation between prior and posterior.

Each criterion arises from its own combination of these three ingredients, and for T-, L-, and C-optimality the Gaussian prior covariance is chosen so that it encodes the criterion's weight matrix. The five criteria are therefore special cases of a single Bayesian-geometric template. The main contribution is presenting all five classical criteria under this template; the T-, L-, C-, and A-optimality derivations appear novel to the best of my knowledge.

A relaxation of the standard setup is then developed. The assumption that the designer knows $\theta_0$ is relaxed by distinguishing the \emph{designer's prior}, used to average over plausible values of $\theta_0$ at design time, from the \emph{modeler's prior}, used for posterior inference; when both priors coincide and the dissimilarity is the Kullback--Leibler (KL) divergence of the posterior relative to the prior, this recovers the classical expected information gain (EIG) objective.

The rest of the paper is organized as follows. Section~\ref{sec:related-works} situates the contribution within the Bayesian experimental design literature. Section~\ref{sec:setup} states the asymptotic setup and the BvM limit used throughout. Section~\ref{sec:framework} introduces the expected prior--posterior dissimilarity objective and its prior-predictive approximation, the classical EIG. Sections~\ref{sec:d-opt}--\ref{sec:a-opt} derive five classical criteria as (approximate) leading-order limits of expected dissimilarities: D-optimality from $\KL(\text{posterior} \parallel \text{prior})$ (Section~\ref{sec:d-opt}), T- and weighted-trace optimality from $\KL(\text{prior} \parallel \text{posterior})$ (Section~\ref{sec:t-opt-revkl}), and L-, C- and A-optimality from a first-order approximation of the squared 2-Wasserstein distance (Section~\ref{sec:a-opt}). Section~\ref{sec:robust} shows that E-optimality arises as a min-max robustification of C-optimality against uncertainty in the modeler's prior direction. Section~\ref{sec:designer} relaxes the assumption that the designer knows $\theta_0$ by introducing a designer's prior, recovering the classical EIG objective when the designer's and modeler's priors coincide. Section~\ref{sec:conclusion} concludes.

\section{Related works}
\label{sec:related-works}

Selecting experiments by maximizing an expected utility, such as the dissimilarity between prior and posterior, has a long history in Bayesian optimal design. \citet{lindley1956} formulated experimental design as maximizing an expected utility, with a natural choice being the expected reduction in uncertainty about $\theta$, expressed via entropy or equivalently the KL divergence between prior and posterior. This information-gain perspective became dominant in Bayesian optimal design, and the resulting expected utility is the classical EIG (formally defined in Eq.~\ref{eq:eig-mi} below); under a Gaussian prior, maximizing the EIG recovers D-optimality in the large-$N$ limit \citep{stone1959,chaloner1995}, as re-derived in Section~\ref{sec:d-opt}. A parallel thread, also rooted in Bayesian decision theory, measures discrepancy through squared Euclidean distance in parameter space, giving rise to the negative squared distance (NSD) criterion,
\begin{align*}
U_{\text{NSD}}(\delta) &= -\E_{\theta, D \mid \delta}\bigl[\|\theta - \mu_1(D, \delta)\|^2\bigr] \\
&= -\E_{D \mid \delta}\bigl[\tr(\Var(\theta \mid D, \delta))\bigr],
\end{align*}
where $\mu_1(D, \delta) = \E[\theta \mid D, \delta]$ denotes the posterior mean, $\theta$ is drawn from the prior $p(\theta)$, and $D$ from the model $p(D \mid \theta, \delta)$ rather than from the true data-generating distribution $p(D \mid \theta_0, \delta)$. Under the BvM limit, $\Var(\theta \mid D, \delta) \approx N^{-1} F(\theta_0, \delta)^{-1}$, so maximizing the NSD asymptotically recovers A-optimality. Both threads were unified under the general expected-utility framework reviewed by \citet{chaloner1995}, where a design is chosen to maximize an expected posterior improvement, whether information-theoretic or geometric. The present paper continues in this tradition: it offers an alternative, transport-based route to A-optimality via an approximate 2-Wasserstein distance (Section~\ref{sec:a-opt}) and reads T-, L-, and C-optimality as further leading-order limits of approximate expected prior--posterior dissimilarities.

A note on terminology: the label ``T-optimality'' is used in this paper for the trace of the FIM, in line with one common convention. A different criterion, also called T-optimality, was introduced by \citet{atkinson1975a} and \citet{atkinson1975b} for discriminating between rival models; the two criteria share a name but address different problems and should not be confused.

The Bayesian experimental design literature is large and has explored information-theoretic, expected-utility, and transport-based criteria from many angles \citep{lindley1956,degroot1962,chaloner1995,ryan2016,rainforth2024,helin2025}, with modern computational developments built on simulation-based and variational EIG estimation \citep{huan2013,foster2019}.

The D- and A-criteria are themselves unified on the classical side by the $\Phi_p$-criterion family of \citet{kiefer1974}, for $p \in (-\infty, \infty)$, which recovers D-optimality in the limit as $p \to 0$, T-optimality at $p = 1$, A-optimality at $p = -1$, and E-optimality in the limit as $p \to -\infty$. In the current paper, a different unification is offered on the Bayesian-geometric side, where D- and A-optimality arise from different dissimilarities rather than from a single family of criteria. Whether additional dissimilarities can recover the intermediate $\Phi_p$ criteria for other values of $p$ remains an open question.

One recent closely related work is \citet{helin2025}, who, restricting attention to the linear Gaussian model, show that minimizing the 2-Wasserstein distance between the prior and posterior is equivalent to minimizing the sum of the square roots of the eigenvalues of the posterior covariance matrix, equivalent to $p = -1/2$ in the $\Phi_p$ family. That derivation is not approximate (unlike the one in Section~\ref{sec:a-opt}) but it is confined to the linear Gaussian setting, does not invoke the BvM theorem to extend to nonlinear models in the large-$N$ limit, and does not account for the data-dependent posterior mean that arises once $D$ is treated as random under the true data-generating distribution.

In a different direction, the robust-design and simulation-based design literature \citep{pronzato1985,muller1999} has developed expected-utility formulations under a designer's prior over $\theta_0$, with \citet{pronzato1985} in particular treating expected functionals of the FIM as design objectives. Section~\ref{sec:designer} places the dissimilarity-based objectives inside this standard framework.

\section{Asymptotic setup}
\label{sec:setup}

Assume the parameters $\theta$ of the model are identifiable and that the standard regularity conditions hold under which $F(\theta_0, \delta)$ is positive definite and the score function is well behaved. Under these conditions, the maximum likelihood estimator $\hat{\theta}_N(D, \delta)$ is asymptotically Gaussian,
\begin{equation}
\label{eq:mle-cov}
\hat{\theta}_N(D, \delta) \,\dot{\sim}\, \N\!\bigl(\theta_0, N^{-1} F(\theta_0, \delta)^{-1}\bigr).
\end{equation}
The BvM theorem \citep{bernstein1917,vonmises1931,lecam1953,vandervaart1998,fischer2025} implies that for any prior $p(\theta)$ that is positive and continuous in a neighborhood of $\theta_0$, the posterior is asymptotically Gaussian, centered at the maximum likelihood estimator, for a large sample size $N$,
\begin{equation}
\label{eq:bvm}
p(\theta \mid D, \delta) \approx \N\!\bigl(\hat{\theta}_N(D, \delta), N^{-1} F(\theta_0, \delta)^{-1}\bigr).
\end{equation}
The asymptotic posterior \emph{covariance} is independent of $D$ to leading order, but the asymptotic posterior \emph{mean} is not.

\section{Expected dissimilarity between prior and posterior}
\label{sec:framework}

Recall from Section~\ref{sec:introduction} the expected prior--posterior dissimilarity
\[
U(\delta) = \E_{D \mid \theta_0, \delta}\bigl[d(p(\theta), p(\theta \mid D, \delta))\bigr],
\]
where the outer expectation is taken under the true data-generating distribution $p(D \mid \theta_0, \delta)$ and not under the model's prior predictive $p(D \mid \delta)$. An optimal design is one that maximizes this expected dissimilarity for some choice of $d$. As noted in Section~\ref{sec:introduction}, the same construction applies to any derived quantity $q = g(\theta)$ by replacing $p(\theta)$ and $p(\theta \mid D, \delta)$ with the prior predictive $p(q)$ and posterior predictive $p(q \mid D, \delta)$; we develop the parameter-space case below, but the derivations carry over to any such $q$.

Since we do not know the true parameters $\theta_0$ at design time, we cannot compute $U(\delta)$ exactly. In the Bayesian optimal design literature, it is common to instead compute $U(\delta)$ under the model's own prior predictive distribution $p(D \mid \delta)$. This is a convenient approximation that allows us to compute $U(\delta)$ without knowing $\theta_0$, but it is important to keep in mind that this approximation presupposes that the model is correctly specified and that the prior is not too misspecified. We denote the approximate expected dissimilarity under the model's own prior predictive distribution as
\[
U_{\text{prior}}(\delta) = \E_{D \mid \delta}\bigl[d(p(\theta), p(\theta \mid D, \delta))\bigr],
\]
where $D$ is assumed to be sampled from the prior predictive distribution $p(D \mid \delta)$. If we further define $d(p_0, p_1) = \KL(p_1 \parallel p_0)$ as the KL divergence \citep{kullback1951} of the posterior relative to the prior, $U_{\text{prior}}(\delta)$ becomes the classical EIG \citep{lindley1956}:
\begin{align}
\label{eq:eig-mi}
U_{\text{prior}}(\delta)
&= \E_{D \mid \delta}\bigl[\KL(p(\theta \mid D, \delta) \parallel p(\theta))\bigr] \notag\\
&= H(\theta) - \E_{D \mid \delta}\bigl[H(\theta \mid D, \delta)\bigr] \geq 0,
\end{align}
where $H(\theta) = -\E_{\theta}[\log p(\theta)]$ is the prior entropy, and $H(\theta \mid D, \delta) = -\E_{\theta \mid D, \delta}[\log p(\theta \mid D, \delta)]$ is the posterior entropy. Reading entropy as a measure of uncertainty,\footnote{%
Entropy is not the only notion of uncertainty in a distribution; variance is another common one, and the two need not move together. Consider the bimodal mixture $p(\theta) = 0.5\,\N(-\alpha, \epsilon^2) + 0.5\,\N(\alpha, \epsilon^2)$, with $\alpha > 1$ and $\epsilon \ll 1$. Provided the modes are well separated ($\alpha \gg \epsilon$), the entropy is approximately $\log 2 + \tfrac{1}{2}\log(2\pi e \epsilon^2)$, an increasing function of $\epsilon$ that does not depend on $\alpha$. The variance, by contrast, is exactly $\alpha^2 + \epsilon^2$, increasing in both $\alpha$ and $\epsilon$. Decreasing $\epsilon$ while increasing $\alpha$ therefore reduces the entropy and increases the variance simultaneously: the two notions of uncertainty move in opposite directions. The entropy interpretation of the EIG as a reduction in uncertainty therefore becomes more subtle when the prior or posterior is non-Gaussian, which is typical when the model is nonlinear or has structural non-identifiabilities, both of which can generate multi-modal posteriors.
} the EIG is the expected reduction in uncertainty about $\theta$ from observing $D$ under design $\delta$, and an optimal design is one that maximizes this expected reduction. The non-negativity in Eq.~\ref{eq:eig-mi} is a direct consequence of the non-negativity of the KL divergence: the \emph{expected} posterior entropy is no larger than the prior entropy under \emph{any} design; a poor $\delta$ can drive the gap to zero, but never below.

This guarantee can seem counter-intuitive when the prior is over-confident and badly placed: one would expect data generated under any plausible $\delta$ to be inconsistent with the bad prior, making the posterior more diffuse than the prior, i.e.\ $H(\theta \mid D, \delta) > H(\theta)$. This apparent contradiction is resolved by noting that:
\begin{enumerate}
\item The expectation is taken under the model's \emph{own} prior predictive distribution $p(D \mid \delta)$, so the non-negativity presupposes that both the prior and the likelihood are correctly specified. Under prior misspecification, the actual data are drawn from some $p(D \mid \theta_0, \delta) \neq p(D \mid \delta)$, so $U_{\text{prior}}(\delta)$ would be averaging dissimilarities against a fictitious data distribution, decoupling it from how much an experiment would actually teach us about $\theta_0$.
\item Even under correct specification, the non-negativity is a statement about the \emph{expectation} over $D$, not about any particular realization. Individual $D$ can produce a posterior whose entropy \emph{exceeds} the prior's; only the mean over the prior predictive distribution is guaranteed not to increase.
\end{enumerate}

\section{D-optimality from $\KL(\text{posterior} \parallel \text{prior})$}
\label{sec:d-opt}

Let the dissimilarity function be the forward KL divergence of the posterior relative to the prior,
\[
d(p_0, p_1) = \KL(p_1 \parallel p_0).
\]
The optimal design in this case is:
\begin{align*}
\delta^* &= \argmax_{\delta \in \D} \E_{D \mid \theta_0, \delta}\bigl[d(p(\theta), p(\theta \mid D, \delta))\bigr] \\
&= \argmax_{\delta \in \D} \E_{D \mid \theta_0, \delta}\bigl[\KL(p(\theta \mid D, \delta) \parallel p(\theta))\bigr].
\end{align*}

Using the asymptotic posterior Eq.~\ref{eq:bvm}, the divergence is approximately
\begin{multline*}
\KL(p(\theta \mid D, \delta) \parallel p(\theta)) \\
\approx \KL\!\Bigl(\N(\hat{\theta}_N(D, \delta), N^{-1} F(\theta_0, \delta)^{-1})
\parallel p(\theta)\Bigr).
\end{multline*}

Take a Gaussian prior $p(\theta) = p_0 = \N(\mu_0, \Sigma_0)$ with arbitrary mean $\mu_0$ and arbitrary positive-definite covariance $\Sigma_0$. Let the asymptotic posterior be $p_1 = \N(\mu_1, \Sigma_1)$ with $\mu_1 = \hat{\theta}_N(D, \delta)$ and $\Sigma_1 = N^{-1} F(\theta_0, \delta)^{-1}$. The KL divergence in this case is:
\begin{multline*}
\KL(p_1 \parallel p_0) = \tfrac{1}{2}\Bigl(
\log\tfrac{\det \Sigma_0}{\det \Sigma_1} - k
+ \tr(\Sigma_0^{-1}\Sigma_1) \\
+ (\mu_1 - \mu_0)^T \Sigma_0^{-1} (\mu_1 - \mu_0)
\Bigr).
\end{multline*}

Substituting the asymptotic posterior parameters, the trace term is
\[
\tr(\Sigma_0^{-1}\Sigma_1) = N^{-1} \tr(\Sigma_0^{-1} F(\theta_0, \delta)^{-1}),
\]
the log-determinant term, using $\det(cA) = c^k \det A$ for a $k \times k$ matrix $A$, is
\begin{align*}
\log\frac{\det \Sigma_0}{\det \Sigma_1}
&= \log \det \Sigma_0 + \log \det(N F(\theta_0, \delta)) \notag\\
&= \log \det \Sigma_0 + k \log N + \log \det F(\theta_0, \delta),
\end{align*}
and the mean-difference term is the only one that depends on $D$. Taking the outer expectation and using Eq.~\ref{eq:mle-cov},
\begin{multline*}
\E_{D \mid \theta_0, \delta}\bigl[(\hat{\theta}_N - \mu_0)^T \Sigma_0^{-1} (\hat{\theta}_N - \mu_0)\bigr] \\
= (\theta_0 - \mu_0)^T \Sigma_0^{-1} (\theta_0 - \mu_0)
+ \tfrac{1}{N} \tr(\Sigma_0^{-1} F(\theta_0, \delta)^{-1}),
\end{multline*}
of which the first piece does not depend on $\delta$ and the second contributes a design-dependent $O(N^{-1})$ term. Collecting all design-dependent contributions of $\E_{D \mid \theta_0, \delta}[\KL(p_1 \parallel p_0)]$,
\begin{align*}
U(\delta) &= \underbrace{\tfrac{1}{2}\log \det F(\theta_0, \delta)}_{O(1)} \notag\\
&\quad + \underbrace{\tfrac{1}{N}\tr(\Sigma_0^{-1} F(\theta_0, \delta)^{-1})}_{O(N^{-1})}
+ \text{const},
\end{align*}
where ``const'' collects all terms that do not depend on $\delta$. The log-determinant dominates for large $N$, so
\[
U(\delta) \approx \tfrac{1}{2}\log \det F(\theta_0, \delta) + \text{const},
\]
and
\begin{align*}
\delta^* &= \argmax_{\delta \in \D} U(\delta) \notag\\
&\approx \argmax_{\delta \in \D} \log \det F(\theta_0, \delta),
\end{align*}
which is D-optimality. This is a classical result \citep{stone1959,chaloner1995}. One can alternatively write D-optimality as:
\begin{align*}
\delta^*_D &= \argmax_{\delta \in \D} \bigl[\tfrac{1}{2}\log \det F(\theta_0, \delta)\bigr] \notag\\
&\approx \argmax_{\delta \in \D} \bigl[U(\delta) - \tfrac{1}{N}\tr(\Sigma_0^{-1} F(\theta_0, \delta)^{-1})\bigr],
\end{align*}
exposing the negative trace term as a sub-leading $O(N^{-1})$ correction that vanishes as $N \to \infty$ and so does not enter the asymptotic criterion.

\section{T- and weighted-trace optimality from $\KL(\text{prior} \parallel \text{posterior})$}
\label{sec:t-opt-revkl}

This section recovers T-optimality and its weighted-trace generalization. Throughout this paper, T-optimality means the trace-of-FIM criterion
\[
\delta^*_T = \argmax_{\delta \in \D} \tr(F(\theta_0, \delta)),
\]
which is one of the two distinct criteria sharing the name in the literature; it is \emph{not} the model-discrimination criterion of \citet{atkinson1975a}. The weighted-trace generalization is
\[
\delta^*_{\text{wT}} = \argmax_{\delta \in \D} \tr(L F(\theta_0, \delta) L^T)
\]
for a weight matrix $L$, and reduces to T-optimality when $L = \I$.

KL divergence is asymmetric, and reversing its arguments isolates a different leading order. Let the dissimilarity function be the KL divergence of the prior relative to the posterior, scaled by $1/N$,
\[
d(p_0, p_1) = \tfrac{1}{N}\,\KL(p_0 \parallel p_1).
\]
Using the same Gaussian prior and asymptotic posterior as above, the un-scaled divergence is
\begin{align*}
&\KL(p(\theta) \parallel p(\theta \mid D, \delta)) \\
&\quad = \tfrac{1}{2}\Bigl(
\log\tfrac{\det \Sigma_1}{\det \Sigma_0} - k
+ \tr(\Sigma_1^{-1}\Sigma_0) \\
&\qquad + (\mu_0 - \hat{\theta}_N)^T \Sigma_1^{-1}
(\mu_0 - \hat{\theta}_N)
\Bigr).
\end{align*}

Using the asymptotic posterior parameters, the trace term can be written as
\[
\tr(\Sigma_1^{-1}\Sigma_0) = N \tr(F(\theta_0, \delta)\,\Sigma_0),
\]
the log-determinant term is
\begin{align*}
\log\frac{\det \Sigma_1}{\det \Sigma_0}
&= \log \det(N^{-1} F(\theta_0, \delta)^{-1}) - \log \det \Sigma_0 \notag\\
&= -\log \det \Sigma_0 - k \log N - \log \det F(\theta_0, \delta),
\end{align*}
and the mean-difference term is the only one carrying $D$-dependence. Using the asymptotic distribution of the MLE and taking the outer expectation, using Eq.~\ref{eq:mle-cov} and the identity $a^T M a = \tr(M a a^T)$,
\begin{align*}
&\E_{D \mid \theta_0, \delta}\bigl[(\mu_0 - \hat{\theta}_N)^T \Sigma_1^{-1} (\mu_0 - \hat{\theta}_N)\bigr] \notag\\
&\qquad = N (\mu_0 - \theta_0)^T F(\theta_0, \delta) (\mu_0 - \theta_0) + k \notag\\
&\qquad = N \tr(F(\theta_0, \delta)\,(\mu_0 - \theta_0)(\mu_0 - \theta_0)^T) + k.
\end{align*}
Collecting the design-dependent contributions of $\E_{D \mid \theta_0, \delta}[\KL(p(\theta) \parallel p(\theta \mid D, \delta))]$ and dividing by $N$,
\begin{align*}
U(\delta) &= \E_{D \mid \theta_0, \delta}\bigl[d(p(\theta), p(\theta \mid D, \delta))\bigr] \\
&= \tfrac{1}{2}\tr(F(\theta_0, \delta)\,\Sigma_0) \\
&\quad + \tfrac{1}{2}(\mu_0 - \theta_0)^T F(\theta_0, \delta) (\mu_0 - \theta_0) \\
&\quad - \tfrac{1}{2N}\log \det F(\theta_0, \delta)
+ \text{const},
\end{align*}
where ``const'' collects all terms that do not depend on $\delta$.

The first term, $\tfrac{1}{2}\tr(F(\theta_0, \delta)\,\Sigma_0)$, is a weighted trace of the FIM with the prior covariance $\Sigma_0$ playing the role of the weight matrix. The second term,
\begin{align*}
R(\delta) &= \tfrac{1}{2}(\mu_0 - \theta_0)^T F(\theta_0, \delta) (\mu_0 - \theta_0) \notag\\
&= \tfrac{1}{2}\tr(F(\theta_0, \delta)\,(\mu_0 - \theta_0)(\mu_0 - \theta_0)^T),
\end{align*}
is non-negative because $F(\theta_0, \delta)$ is positive-definite. Dropping the second term yields a lower bound
\[
U(\delta) \geq \tfrac{1}{2}\tr(F(\theta_0, \delta)\,\Sigma_0)
- \tfrac{1}{2N}\log \det F(\theta_0, \delta)
+ \text{const}.
\]
For large $N$, the $O(N^{-1})$ log-determinant term is sub-leading, and maximizing the lower bound reduces to
\[
\delta^* \approx \argmax_{\delta \in \D} \tr(F(\theta_0, \delta)\,\Sigma_0),
\]
a weighted-trace criterion with weight matrix $\Sigma_0$. Taking $\Sigma_0 = \I$ gives plain T-optimality. Thus T-optimality arises as the limit of a lower bound on the reverse-KL objective.

The two terms can also be combined:
\begin{multline*}
\tr(F(\theta_0, \delta)\,\Sigma_0)
+ \tr(F(\theta_0, \delta)\,(\mu_0 - \theta_0)(\mu_0 - \theta_0)^T) \\
= \tr(F(\theta_0, \delta)\,M),
\end{multline*}
where $M = \Sigma_0 + (\mu_0 - \theta_0)(\mu_0 - \theta_0)^T$. Since $M$ is positive-definite, it admits a factorization $M = L^T L$, yielding
\[
\tr(F(\theta_0, \delta)\,M) = \tr(L F(\theta_0, \delta) L^T).
\]
The full expected reverse-KL therefore has $\tfrac{1}{2}\tr(L F(\theta_0, \delta) L^T)$ as its leading-order term, with the weight matrix $L$ encoding both the prior covariance and the mean-offset. Notably, when the prior's covariance is sufficiently diffuse
\[
\Sigma_0 + (\mu_0 - \theta_0)(\mu_0 - \theta_0)^T \approx \Sigma_0,
\]
resulting in the lower bound being approximately tight and the leading-order term being approximately $\tfrac{1}{2}\tr(F(\theta_0, \delta)\,\Sigma_0)$, so the weighted-trace criterion is recovered even when the prior mean is not centered at the truth, provided the prior covariance is sufficiently large.

\section{L-, C- and A-optimality from 2-Wasserstein distance}
\label{sec:a-opt}

Let the dissimilarity function be the squared 2-Wasserstein distance \citep{villani2009},
\[
d(p_0, p_1) = N \cdot W_2^2(p_0, p_1).
\]
The Wasserstein distance is computed between the Gaussian prior $p_0 = \N(\mu_0, \Sigma_0)$ and the BvM Gaussian limit of the posterior, $p_1 = \N(\hat{\theta}_N(D, \delta), N^{-1} F(\theta_0, \delta)^{-1})$. For two Gaussians, there is a closed form \citep{dowson1982,olkin1982,givens1984}:
\begin{multline*}
W_2^2(p_0, p_1) = \|\mu_0 - \mu_1\|^2 \\
+ \tr\!\left(\Sigma_0 + \Sigma_1 - 2(\Sigma_0^{1/2}\Sigma_1\Sigma_0^{1/2})^{1/2}\right).
\end{multline*}

To linearize the matrix-square-root term, recall that the Sylvester equation $AY + YB = Z$ admits a unique solution $Y$ if and only if the spectra of $A$ and $-B$ are disjoint, $\sigma(A) \cap \sigma(-B) = \emptyset$ \citep[Theorem VII.2.1]{bhatia1997}. For the Lyapunov form $AY + YA = Z$ with positive-definite $A$, one has $\sigma(A) \subset (0, \infty)$ and $\sigma(-A) \subset (-\infty, 0)$, so this disjointness condition holds automatically and the inverse Lyapunov operator $L_A^{-1}[Z] = Y$ is well-defined. The Fr\'echet derivative of the matrix square root at a positive-definite $X_1$ is precisely $L_{X_1^{1/2}}^{-1}$ \citep{higham2008}. Hence
\[
X^{1/2} = X_1^{1/2} + L_{X_1^{1/2}}^{-1}[X - X_1] + O(\|X - X_1\|^2).
\]

Substituting $X_1 = \epsilon^2 \I$ for some $\epsilon > 0$ gives $X_1^{1/2} = \epsilon \I$ and the Sylvester equation $\epsilon Y + Y \epsilon = Z$ has solution $Y = Z/(2\epsilon)$, so $L_{X_1^{1/2}}^{-1}[Z] = Z/(2\epsilon)$ and the square root term is
\begin{multline*}
(\Sigma_0^{1/2}\Sigma_1\Sigma_0^{1/2})^{1/2}
= \epsilon \I + \tfrac{1}{2\epsilon}(\Sigma_0^{1/2}\Sigma_1\Sigma_0^{1/2} - \epsilon^2 \I) \\
+ O(\|\Sigma_0^{1/2}\Sigma_1\Sigma_0^{1/2} - \epsilon^2 \I\|^2).
\end{multline*}
For $\epsilon$ small enough but positive, this linearization is accurate at large $N$: the posterior covariance has norm $\|\Sigma_1\| = O(N^{-1})$, so the matrix $\Sigma_0^{1/2}\Sigma_1\Sigma_0^{1/2}$ being expanded vanishes with $N$ and $\epsilon$ can be picked small enough that $\|\Sigma_0^{1/2}\Sigma_1\Sigma_0^{1/2} - \epsilon^2 \I\|$ is small.

Dropping the higher order term, the Wasserstein distance is approximately
\begin{align*}
W_2^2(p_0, p_1)
&\approx \|\mu_0 - \mu_1\|^2
+ \tr\!\left(\Sigma_0 + \Sigma_1 - 2\epsilon \I\right. \notag\\
&\qquad \left.{}- \tfrac{1}{\epsilon}(\Sigma_0^{1/2}\Sigma_1\Sigma_0^{1/2} - \epsilon^2 \I)\right) \notag\\
&= \|\mu_0 - \mu_1\|^2 + \tr(\Sigma_0 + \Sigma_1 - \epsilon \I) \notag\\
&\quad - \tfrac{1}{\epsilon}\tr(\Sigma_0^{1/2}\Sigma_1\Sigma_0^{1/2}).
\end{align*}
Using the linear and cyclic properties of the trace, the last term is
\[
\tr(\Sigma_0^{1/2}\Sigma_1\Sigma_0^{1/2}) = \tr(\Sigma_1\Sigma_0) = \tr(\Sigma_0\Sigma_1),
\]
and the Wasserstein distance is approximately
\begin{align*}
W_2^2(p_0, p_1)
&\approx \|\mu_0 - \mu_1\|^2
+ \tr(\Sigma_0 + \Sigma_1 - \epsilon \I - \tfrac{1}{\epsilon}\Sigma_0\Sigma_1) \notag\\
&= \|\mu_0 - \mu_1\|^2 + \tr(\Sigma_0) - \epsilon k \notag\\
&\quad + \tr\!\left((\I - \tfrac{1}{\epsilon}\Sigma_0)\Sigma_1\right).
\end{align*}

Substituting the asymptotic posterior parameters $\mu_1 = \hat{\theta}_N(D, \delta)$ and $\Sigma_1 = N^{-1} F(\theta_0, \delta)^{-1}$,
\begin{multline*}
W_2^2(p(\theta), p(\theta \mid D, \delta))
\approx \|\mu_0 - \hat{\theta}_N(D, \delta)\|^2 \\
+ \tr(\Sigma_0) - \epsilon k
+ N^{-1} \tr\!\left((\I - \tfrac{1}{\epsilon}\Sigma_0) F(\theta_0, \delta)^{-1}\right).
\end{multline*}
The first term is the only one with $D$-dependence. Taking the outer expectation over $D$ and using Eq.~\ref{eq:mle-cov},
\begin{align*}
&\E_{D \mid \theta_0, \delta}\bigl[\|\mu_0 - \hat{\theta}_N(D, \delta)\|^2\bigr] \notag\\
&\quad = \|\mu_0 - \theta_0\|^2
+ N^{-1} \tr(F(\theta_0, \delta)^{-1}),
\end{align*}
of which the first piece does not depend on $\delta$ and the second is a design-dependent $O(N^{-1})$ contribution. Collecting all design-dependent contributions,
\begin{multline*}
\E_{D \mid \theta_0, \delta}\bigl[W_2^2(p(\theta), p(\theta \mid D, \delta))\bigr] \\
\approx N^{-1} \tr\!\left((2\I - \tfrac{1}{\epsilon}\Sigma_0) F(\theta_0, \delta)^{-1}\right) + \text{const}.
\end{multline*}
where ``const'' collects all terms that do not depend on $\delta$.

Let $\Sigma_0 = 2\epsilon \I + \epsilon B$ for some arbitrary positive semi-definite matrix $B$. This parameterization does not lose generality because for a given prior covariance $\Sigma_0$ and $\epsilon$, we can always solve for $B$ so that $\Sigma_0 = 2\epsilon \I + \epsilon B$ holds. Additionally, the result below is agnostic to the scale of $B$ as long as it does not depend on $N$. The above means that $2\I - \tfrac{1}{\epsilon}\Sigma_0 = -B$. The expected Wasserstein distance is then approximately
\begin{multline*}
\E_{D \mid \theta_0, \delta}\bigl[W_2^2(p(\theta), p(\theta \mid D, \delta))\bigr] \\
\approx -N^{-1} \tr(B F(\theta_0, \delta)^{-1}) + \text{const}.
\end{multline*}
Since $B$ is positive semi-definite, we can always decompose it as $B = L^T L$ for some matrix $L$ (does not have to be triangular). By the cyclic property of the trace, the expected Wasserstein distance is approximately
\begin{multline*}
\E_{D \mid \theta_0, \delta}\bigl[W_2^2(p(\theta), p(\theta \mid D, \delta))\bigr] \\
\approx -N^{-1} \tr(L F(\theta_0, \delta)^{-1} L^T) + \text{const}.
\end{multline*}
The optimal design is then
\begin{align*}
\delta^* &= \argmax_{\delta \in \D} U(\delta) \notag\\
&= \argmax_{\delta \in \D} \E_{D \mid \theta_0, \delta}\bigl[N \cdot W_2^2(p(\theta), p(\theta \mid D, \delta))\bigr] \notag\\
&\approx \argmax_{\delta \in \D} -\tr(L F(\theta_0, \delta)^{-1} L^T) \notag\\
&= \argmin_{\delta \in \D} \tr(L F(\theta_0, \delta)^{-1} L^T).
\end{align*}
This is known as L-optimality \citep{pukelsheim1993}. When $L$ is an arbitrary positive scalar $\alpha$ times the identity matrix, this reduces to A-optimality:
\begin{align*}
\delta^*_A &= \argmin_{\delta \in \D} \alpha^2 \tr(F(\theta_0, \delta)^{-1}) \notag\\
&= \argmin_{\delta \in \D} \tr(F(\theta_0, \delta)^{-1}).
\end{align*}
If $B$ is instead rank-1
\[
B = \alpha^2 \cdot c \cdot c^T,
\]
where $\alpha > 0$ is an arbitrary positive scalar and $c$ is a unit vector, then setting $L = \alpha c^T$ reduces to C-optimality:
\begin{align*}
\delta^*_C &= \argmin_{\delta \in \D} \alpha^2 c^T F(\theta_0, \delta)^{-1} c \notag\\
&= \argmin_{\delta \in \D} c^T F(\theta_0, \delta)^{-1} c.
\end{align*}
Therefore, L-, A- and C-optimality arise as the leading-order term of the expected approximate 2-Wasserstein distance between the prior and the BvM Gaussian posterior, under a Gaussian prior.

C-optimality also admits a predictive-distribution interpretation: applying the derived-quantity construction of Section~\ref{sec:framework} to the scalar $q = c^T \theta$ (the linear functional originally targeted by \citet{elfving1952}), it can be read as the optimal design for maximizing the expected dissimilarity between the prior predictive $p(q)$ and posterior predictive $p(q \mid D, \delta)$ of that quantity, with the same $c^T F(\theta_0, \delta)^{-1} c$ objective emerging as the leading-order term.

\section{Robustness to the modeler's prior: E-optimality as robust C-optimality}
\label{sec:robust}

The L-/C-/A-optimality derivation in Section~\ref{sec:a-opt} parametrizes the prior covariance as $\Sigma_0 = 2\epsilon \I + \epsilon B$, where $B = L^T L$ encodes the criterion's weight matrix. C-optimality arises when $B$ is rank-1, so $B = \alpha^2 \cdot c \cdot c^T$ for some unit vector $c$ and positive scalar $\alpha$, yielding the objective $c^T F(\theta_0, \delta)^{-1} c$. The direction $c$ and scalar $\alpha$ are determined by the modeler's prior.

A natural robust-design question arises when the designer does not know, or does not fully trust, the modeler's choice of direction $c$. Suppose the modeler believes that only a single linear combination $c^T \theta$ matters, but the designer is uncertain about which combination. The designer may then seek a design that performs well in the worst case over plausible directions. If we assume that any direction $c$ with $\|c\| = 1$ is possible, the robustified C-optimality objective becomes
\[
\delta^* = \argmin_{\delta \in \D} \max_{\|c\|=1} c^T F(\theta_0, \delta)^{-1} c.
\]
By the Rayleigh quotient theorem, the inner maximum is the largest eigenvalue of $F(\theta_0, \delta)^{-1}$:
\begin{align*}
\max_{\|c\|=1} c^T F(\theta_0, \delta)^{-1} c
&= \lambda_{\max}(F(\theta_0, \delta)^{-1}) \notag\\
&= \frac{1}{\lambda_{\min}(F(\theta_0, \delta))}.
\end{align*}
The robust objective is therefore
\[
\delta^* = \argmax_{\delta \in \D} \lambda_{\min}(F(\theta_0, \delta)),
\]
which is exactly E-optimality.

This derivation places E-optimality not as an independent Bayesian-geometric criterion, but as the robustification of C-optimality against uncertainty in the modeler's prior direction. Where C-optimality assumes the modeler knows exactly which linear functional $c^T \theta$ is of interest, E-optimality assumes the designer guards against the worst-case choice of $c$, a connection known in the classical literature \citep{kiefer1974}.

\section{Designer's prior}
\label{sec:designer}

The objective $U(\delta)$ in Section~\ref{sec:framework} was written as an expectation under $p(D \mid \theta_0, \delta)$, presupposing that the designer knows the true parameters $\theta_0$. In practice, the designer does not know $\theta_0$ but may have a prior distribution over it. Notably, this prior does not have to coincide with the modeler's prior $p(\theta)$, which is used in Bayes' rule to form the posterior \emph{after} observing $D$. There are therefore two distinct prior distributions in play:
\begin{enumerate}
\item The \emph{modeler's prior} $p(\theta)$, which encodes the modeler's beliefs about $\theta$ before observing $D$, and is used in Bayes' rule to form the posterior.
\item The \emph{designer's prior} $p_0(\theta)$, which encodes the designer's beliefs about plausible values of the unknown $\theta_0$ at design time, and is used to take an outer expectation over $\theta_0$ in the design objective.
\end{enumerate}
Consider the following thought experiment to clarify the distinction. A teacher is designing a synthetic data experiment to teach a student about Bayesian inference. The teacher chooses a model and a set of parameters $\theta_0$ from which data will be generated. The teacher knows the student always picks some Gaussian prior $p(\theta) = \N(\mu_0, \Sigma_0)$. Using the methods in this paper, the teacher can design an experiment optimal for the student to learn about $\theta_0$ under that prior:
\[
\delta^* = \argmax_{\delta \in \D} \E_{D \mid \theta_0, \delta}\bigl[d(p(\theta), p(\theta \mid D, \delta))\bigr].
\]
In this example, the teacher is the designer, the student is the modeler, $\theta_0$ is the true parameters from which data will be generated, and $p(\theta)$ is the modeler's prior. The teacher knows $\theta_0$ but the student does not, and the teacher designs an experiment to maximize the expected learning of the student about $\theta_0$ under the student's prior. The asymptotic derivations of Sections~\ref{sec:d-opt}--\ref{sec:a-opt} all sit inside this asymmetric setup. \citet{etzioni1993} refer to this as designing an experiment for another's analysis.

In practice, the designer does not know the true parameters but may have a prior distribution over them \citep{pronzato1985,chaloner1995}. Denote this by $p_0(\theta)$ and its prior predictive by $p_0(D \mid \delta) = \int p_0(\theta)\,p(D \mid \theta, \delta)\,\mathrm{d}\theta$. The designer then chooses
\begin{align}
\delta^* &= \argmax_{\delta \in \D} \E_{p_0(D \mid \delta)}\bigl[d(p(\theta), p(\theta \mid D, \delta))\bigr] \notag\\
&= \argmax_{\delta \in \D} \E_{\theta_0 \sim p_0(\theta)}\bigl[\E_{D \mid \theta_0, \delta}\bigl[d(p(\theta), \notag\\
&\qquad\qquad p(\theta \mid D, \delta))\bigr]\bigr].
\label{eq:designer}
\end{align}
For each sample of $\theta_0$ from $p_0(\theta)$, that sample is treated as the true parameters in the inner expectation.

The designer's prior $p_0(\theta)$ may or may not coincide with the modeler's prior $p(\theta)$. If the designer chooses to use the same prior that the designer thinks the modeler will use (e.g.\ because they are the same person) and if the dissimilarity is the KL divergence of the posterior relative to the prior, Eq.~\ref{eq:designer} recovers the classical EIG objective \citep{lindley1956,bernardo1979,chaloner1995}:
\begin{align*}
\delta^* &= \argmax_{\delta \in \D} \E_{\theta_0 \sim p(\theta)} \notag\\
&\qquad \bigl[\E_{D \mid \theta_0, \delta}\bigl[\KL(p(\theta \mid D, \delta) \parallel p(\theta))\bigr]\bigr] \\
&= \argmax_{\delta \in \D} \E_{D \mid \delta}\bigl[\KL(p(\theta \mid D, \delta) \parallel p(\theta))\bigr],
\end{align*}
where $\E_{D \mid \delta}$ denotes expectation under the prior predictive $p(D \mid \delta) = \int p(\theta)\,p(D \mid \theta, \delta)\,\mathrm{d}\theta$. The asymptotic derivations in Sections~\ref{sec:d-opt}--\ref{sec:a-opt} produce design-dependent objectives of the form $g(F(\theta_0, \delta))$; under the designer's prior these are averaged to $\E_{\theta_0 \sim p_0(\theta)}[g(F(\theta_0, \delta))]$, recovering the standard ``Bayesian alphabetic'' criteria of \citet{chaloner1995}.

Averaging over $\theta_0$ under a designer's prior to handle this is a standard move in the robust- and Bayesian-design literature \citep{pronzato1985,chaloner1995}, so there is no claim to novelty in this section. It is included for completeness and to clarify how the alphabetic criteria derived above interact with the designer's prior as well as to highlight the distinction between the modeler's and designer's priors.

\section{Conclusion}
\label{sec:conclusion}

I have laid out five classical model-based optimal design criteria (D-, T-, L-, C-, and A-optimality) as the leading-order term, as $N \to \infty$, of five different expected (approximate) prior--posterior dissimilarities under a Gaussian prior and the BvM asymptotic posterior:

\begin{itemize}
\item the expected $\KL(\text{posterior} \parallel \text{prior})$ reproduces the classical D-optimality result of \citet{stone1959} and \citet{chaloner1995};
\item the expected $\KL(\text{prior} \parallel \text{posterior})$ gives T-optimality (the trace-of-FIM criterion) and its weighted-trace generalization;
\item the expected first-order approximation of the squared 2-Wasserstein distance yields L-, C- and A-optimality.
\end{itemize}

While the D-optimality recovery is a classical result, the derivations of T-, L-, C-, and A-optimality from their respective expected dissimilarities are new. These recoveries place the five criteria as special cases of a single Bayesian-geometric template, each obtained by a different setting of the dissimilarity, prior parametrization, and approximation scheme. A sixth classical criterion, E-optimality, is recovered not as an independent expected-dissimilarity limit but as a min-max robustification of C-optimality against uncertainty in the modeler's prior direction.

An extension relaxes the standard assumptions of the asymptotic setup. The assumption that the designer knows the true parameters $\theta_0$ is replaced by a designer's prior $p_0(\theta)$ used to average over plausible $\theta_0$ at design time; this designer's prior is conceptually distinct from the modeler's prior $p(\theta)$ used for posterior inference, and when the two coincide and the dissimilarity is the KL divergence of the posterior relative to the prior, the resulting objective is the classical EIG. While no claim to novelty is made for this extension, it is included for completeness and to clarify how the alphabetic criteria interact with the designer's prior as well as to highlight the distinction between the modeler's and designer's priors.

A few limitations and open directions are worth noting. First, the asymptotic approximations rely on $N$ being large enough for the BvM theorem to apply; at moderate $N$, the prior shape matters and the Gaussian posterior approximation can be less accurate, especially in nonlinear models with multi-modal posteriors. Second, the Bayesian formulation can be more computationally demanding than classical FIM-based design, since evaluating expected dissimilarities often calls for nested Monte Carlo or variational approximation; whether the Bayesian perspective also enables new computational tools is an interesting question. Third, the connection to Kiefer's $\Phi_p$-family suggests asking whether other values of $p$ can be associated with natural dissimilarity functions. Fourth, the paper assumes identifiability and correct model specification, and extensions to locally non-identifiable or misspecified models would involve modified BvM limits, which we leave to future work. Finally, concrete numerical examples, comparing the finite-$N$ behavior of the Bayesian objectives to their asymptotic limits, would be valuable as a complement to the theoretical derivations.

\bibliographystyle{imsart-nameyear}
\bibliography{../od}

\end{document}